\documentclass[aap,MSNbibl,dvips]{arximspdf}

\doi{10.1214/10-AAP698}
\volume{21}
\issue{1}
\pubyear{2011}
\firstpage{332}
\lastpage{350}

\makeatletter

\newtheorem{theorem}{Theorem}[section]

\newtheorem{proposition}[theorem]{Proposition}
\newproclaim{definition}[theorem]{Definition}
\newtheorem{hypothesis}[theorem]{Hypothesis}
\newtheorem{assumption}[theorem]{Assumption}
\newproclaim{example}{Example}
\newproclaim{remark}[theorem]{Remark}
\makeatother

\begin{document}
\begin{frontmatter}

\title{Boundary conditions for the single-factor term structure equation}
\runtitle{The single-factor term structure equation}

\begin{aug}
\author{\fnms{Erik} \snm{Ekstr\"{o}m}\thanksref{t1,t2}\ead[label=e1]{ekstrom@math.uu.se}\corref{}}
\and
\author{\fnms{Johan} \snm{Tysk}\thanksref{t1}\ead[label=e2]{johan.tysk@math.uu.se}}
\thankstext{t1}{Supported by the Swedish Research Council
(VR).}
\thankstext{t2}{Supported by the Swedish National
Graduate School of Mathematics and Computing (FMB).}
\runauthor{E. Ekstr\"{o}m and J. Tysk}
\affiliation{Uppsala University}
\address{Department of Mathematics\\ Uppsala University\\
Box 480, SE-75106 Uppsala\\ Sweden\\
\printead{e1}\\ \phantom{E-mail: }\printead*{e2}} 
\end{aug}

\received{\smonth{11} \syear{2008}}
\revised{\smonth{12} \syear{2009}}

%
\begin{abstract}
We study the term structure equation for single-factor
models that predict nonnegative short rates.
In particular, we show that
the price of a bond or a bond option is the
unique
classical solution to a parabolic differential equation
with a certain boundary behavior for vanishing values of the
short rate.
If the boundary is attainable then this boundary behavior
serves as a boundary condition and guarantees uniqueness of
solutions. On the other hand, if the boundary is nonattainable
then the boundary behavior is not needed to guarantee
uniqueness but it is nevertheless very useful, for instance,
from a numerical perspective.
\end{abstract}

%

\end{frontmatter}

\section{Introduction}\label{sec1}

When calculating prices of different interest
rate derivatives, such as bonds and bond options,
stochastic methods seem to be more commonly used than PDE methods;
compare for instance \cite{BM}. This is in contrast to the case of
stock option pricing where PDE methods are used extensively,
in particular for low dimensional problems.
We believe that one possible
explanation for this phenomenon is that the correspondence between the
risk neutral valuation approach and the pricing equation
(henceforth referred to as the \textit{term structure equation}) with
appropriate boundary conditions
is not fully developed.

For the Black--Scholes
equation, the boundary condition is of Dirichlet type
which corresponds to the underlying asset being
absorbed if reaching zero; compare \cite{JT}.
In contrast, for most interest rate models this is not the case since
the short rate typically would not stay zero if the value zero is reached.
Consequently, it is not clear what boundary conditions
should be specified for the term structure equation.
In fact, the recent monograph \cite{D}
draws attention to this issue in a section
entitled ``The thorny issue of boundary conditions.''
Moreover, in \cite{HLW}, two different solutions to the
term structure equation are presented in the case of the CIR-model.
They are both bounded
and with the same terminal condition, but
naturally they exhibit different boundary behavior for
vanishing interest rates. The authors of that paper take the
view that these solutions represent alternative possible prices.
We on the other hand regard only the solution given by the stochastic
representation as the price, and the purpose of the present paper
is to identify the boundary condition that this stochastic
representation satisfies.

We consider the classical case of a single-factor
model predicting nonnegative values of the short rate.
More precisely, the rate $X(t)$
is modeled directly under the pricing measure as
\begin{eqnarray*}
dX(t)=\beta(X(t),t) \,dt+ \sigma(X(t),t)\, dW,
\end{eqnarray*}
where $W$ is a Brownian motion and $\sigma(0,t)=0$ and $\beta
(0,t)\geq0$.
As indicated above, the option price $u$ corresponding to a
payoff function $g$ is given using risk neutral valuation by
\[
u(x,t)=E_{x,t}\bigl[e^{-\int_t^T X(s)\, ds}g(X(T))\bigr].
\]
Note that if the payoff $g\equiv1$, then bond prices are obtained.
Also note that the set-up covers the case of
bond options.
The corresponding term structure equation is
\[
u_t(x,t)+\tfrac{1}{2}\sigma^2(x,t) u_{xx}(x,t)+\beta(x,t)u_x(x,t)=xu(x,t)
\]
with terminal condition $u(x,T)=g(x)$.
If the price $u$ is twice continuously differentiable up to
and including the boundary $x=0$, then plugging in $x=0$ in
the equation would give the boundary behavior
%
\begin{equation}
\label{bcpp}
u_t(0,t)+\beta(0,t)u_x(0,t)=0.
\end{equation}
Even though there is extensive literature on equations with
degenerating coefficients, compare the classical reference \cite{OR},
the $C^1$-regularity of $u$ at the boundary is not available in the
generality that is needed here.
In fact, one of the solutions in the example in \cite{HLW}
referred to above
is bounded and continuous, but fails to be $C^1$ up to the boundary.
(This solution, however, is not the one given by stochastic
representation.)

In the present paper, sufficient regularity of
the option price $u$ for (\ref{bcpp})
to hold is
established using the Girsanov theorem and
scaling arguments; see Sections~\ref{u_x} and~\ref{u_xx}.
Let us emphasize that (\ref{bcpp}) is the correct
boundary behavior of the option price
regardless if the boundary is hit with positive probability or not.
If the boundary can be reached with positive probability,
then this boundary behavior serves as a boundary condition
for the term structure equation and guarantees uniqueness.
On the other hand, if the boundary is reached
with probability zero, equation (\ref{bcpp}) is not needed
to identify the solution given by the stochastic representation,
and the term ``boundary condition'' is perhaps misleading.
However, it is still valid and certainly useful, for instance,
from a numerical perspective. Indeed, the results of the present
paper have already been implemented numerically in
\cite{ELT}. For simplicity of the exposition,
we will refer to (\ref{bcpp}) as the boundary condition
regardless if the boundary can be reached or not.

In Section~\ref{halfline} the assumptions on the model
and our main result, Theorem~\ref{main}, are presented.
In Sections~\ref{u_x}, \ref{u_xx} and \ref{sahlin}, we establish
regularity properties of the value
function, and use them to prove Theorem~\ref{main}. Finally,
in Section~\ref{realline} we also provide the link between the
stochastic problem and the term structure equation for models
defined on the whole real line.

\section{Assumptions and the main result}
\label{halfline}

When studying the term structure equation on
the positive real axis, we consider models which are
specified so that the short rate automatically stays nonnegative, that
is, there is no need to impose any boundary
behavior of the underlying diffusion process. Throughout
Sections~\ref{halfline}--\ref{sahlin} we work under the
following hypothesis:

\begin{hypothesis}
\label{hyp}
The drift $\beta\in C([0,\infty)\times[0,T])$ is
continuously differentiable in $x$ with bounded derivative,
and $\beta(0,t)\geq0$ for all $t$.
The volatility $\sigma\in C([0,\infty)\times[0,T])$ is such that
$\alpha(x,t):=\frac{1}{2}\sigma^2(x,t)$ is
continuously differentiable in $x$ with a H\"{o}lder continuous
derivative, and
$\sigma(x,t)=0$ if and only if $x=0$.
The functions $\beta$, $\sigma$ and $\alpha_x$ are all
of, at most, linear growth:
%
\begin{equation}
\label{bound}
\vert\beta(x,t)\vert+\vert\sigma(x,t)\vert+
\vert\alpha_x(x,t)\vert\leq C(1+x)
\end{equation}
for all $x$ and $t$. The payoff function $g\dvtx [0,\infty)\to
[0,\infty)$ is continuously differentiable with both $g$ and
$g^\prime$ bounded.
\end{hypothesis}

Let $W$ be a standard Brownian motion on a filtered
probability space
$(\Omega, \mathcal F, (\mathcal F_t)_{t\geq0},P)$. Since
$\alpha$ is continuously differentiable, $\sigma$ is locally
H\"{o}lder (1$/$2) in $x$. It follows that there
exists a unique strong solution $X(t)$ to
%
\begin{equation}
\label{X}
dX(t)=\beta(X(t),t)\, dt+ \sigma(X(t),t)\, dW
\end{equation}
for any initial point $x\geq0$; compare Section~IX.3 in
\cite{RY}. Moreover, it follows from
monotonicity results for stochastic differential equations
with respect to the drift coefficient
(e.g., Theorem~IX.3.7 in \cite{RY}),
that $X$ remains nonnegative at all times. Indeed,
if $\beta$ is replaced with $\beta\wedge0$, then the
corresponding solution to (\ref{X}) is absorbed at zero,
so $X$ is nonnegative if $\beta(0,t)\geq0$.
The option price $u\dvtx [0,\infty)\times[0,T]\to[0,\infty)$
corresponding to a
payoff function $g\dvtx [0,\infty)\to[0,\infty)$ is given by
%
\begin{equation}
\label{u}
u(x,t)=E_{x,t}\bigl[e^{-\int_t^T X(s) \,ds}g(X(T))\bigr],
\end{equation}
where the indices indicate that $X(t)=x$.
As described in the \hyperref[sec1]{Introduction}, the corresponding
term structure equation is given by
%
\begin{equation}
\label{tse}
u_t(x,t)+\tfrac{1}{2}\sigma^2(x,t) u_{xx}(x,t)+\beta(x,t)u_x(x,t)=xu(x,t)
\end{equation}
for $(x,t)\in(0,\infty)\times[0,T)$, with terminal condition
%
\begin{equation}
\label{terminal}
u(x,T)=g(x).
\end{equation}
Moreover, by formally inserting $x=0$ in the equation we
get the boundary condition
%
\begin{equation}
\label{bc}
u_t(0,t)+\beta(0,t)u_x(0,t)=0
\end{equation}
for all $t\in[0,T)$, since $\sigma(0,t)=0$ by assumption.
One of the main efforts in this paper is to show that
the option price $u$ is continuously differentiable up to the
boundary $x=0$, and that it indeed satisfies the
boundary condition (\ref{bc}) in the classical sense.

\begin{definition}
A \textit{classical solution} to the term structure equation is a function
$v\in C([0,\infty)\times[0,T]) \cap C^1([0,\infty)\times[0,T))\cap
C^{2,1}((0,\infty)\times[0,T))$ which satisfies (\ref{tse}), (\ref
{terminal}) and (\ref{bc}).
\end{definition}

Our main result in this article is the following:

\begin{theorem}
\label{main}
In addition to Hypothesis~$\ref{hyp}$, also assume that
Assumption~\ref{ass} below holds.
The option price $u$ as given by (\ref{u}) is then the unique
bounded classical solution to the term structure equation.
\end{theorem}

\begin{example*}
Classical short rate models such as the Cox--Ingersoll--Ross model
%
\begin{equation}
\label{CIR}
dX(t)=\bigl(a-bX(t)\bigr)\, dt +\sigma\sqrt{X(t)}\, dW,
\end{equation}
and the Dothan model
%
\begin{equation}
\label{D}
dX(t)=aX(t) \,dt+\sigma X(t)\, dW ,
\end{equation}
have boundary conditions at $x=0$ that are immediate to write down.
These conditions are
\[
u_t+au_x=0\quad \mbox{and}\quad u_t=0,
\]
respectively. We note that the boundary condition $u_t=0$ for
the Dothan model means that $u$ is constant along the
boundary, that is,
$u(0,t)=g(0)$. This is the same type of boundary condition that appears
for options on stocks
in \cite{JT}, which can be explained by the fact that the Dothan model
is a geometric Brownian motion.
Theorem~\ref{main} also covers the Hull--White model
%
\begin{equation}
\label{HW}
dX(t)=\bigl(a(t)-b(t)X(t)\bigr)\, dt +\sigma(t)\sqrt{X(t)} \,dW
\end{equation}
(which is a time-dependent generalization of the
Cox--Ingersoll--Ross model), and models of, for example, the form
%
\begin{equation}
\label{exam}
dX(t)=\bigl(b-aX(t)\bigr)\, dt +\sigma X^{\gamma}(t)\, dW,\qquad
\gamma\in(1/2,1],
\end{equation}
which also would be natural to consider for bond
pricing.
\end{example*}

\begin{remark*}
 It seems that many of the classical models for the short
rate are proposed
for their analytical tractability. In particular, if the drift
$\beta$ and the diffusion coefficient $\sigma^2$ are affine,
then the model admits an affine term structure.
It is easy to check that known explicit formulas for bond prices and
bond options satisfy the boundary condition (\ref{bc}).
In particular, for models admitting an affine term structure, it is a
consequence of the associated Riccati equations (see \cite{B},
equation~22.25)
that these boundary conditions are fulfilled.
\end{remark*}

\begin{remark*}
 The assumption that $g$ is continuously
differentiable is
satisfied for bonds, but not in general for bond options.
However, using the Markov property, Theorem~\ref{main}
readily extends to bounded
Lipschitz payoffs provided one can show that the corresponding
option price $x\mapsto u(x,T-\varepsilon)$ is continuously
differentiable on $(0,\infty)$ for any $\varepsilon>0$.
The regularizing effect of parabolic equations guarantees
continuous differentiability on $(0,\infty)$, so
the main difficulty is to show that $x\mapsto u_x(x,T-\varepsilon)$ is
continuous also at 0. If the model is convexity preserving,
this is easily done in certain cases
including, for example, call options written on
bond prices. (Note that the corresponding payoff function
$g$ is bounded since bond prices are bounded.)
For details on which
short rate models are convexity preserving, see \cite{ET}.
To our knowledge, all models used in practice belong to this
class.

One should note that the differentiability of the option price up to
the boundary $x=0$ is not valid
without some Lipschitz bound of $g$ at 0.
To see this, consider the contract function
$g(x)=e^{-2\sqrt{x}}$. Then, with $\beta(x,t)=\frac{1}{2}x$
and $\sigma(x,t)=\sqrt2x$, it is straightforward using
the It\^o formula to show that the process
\[
Y(s)=e^{-\int_t^sX(r)\, dr}g(X(s))
\]
is a martingale. Consequently,
the option price $u$ is given by $u(x,t)=g(x)$ for all $t$,
which fails to be a classical solution to the term structure
equation since it is not differentiable at $x=0$. One might argue,
though, that the boundary condition $u_t(0,t)+
\beta(0,t)u_x(0,t)=0$ is satisfied in a weak sense.
\end{remark*}

The proof of Theorem~\ref{main} is carried out in several steps.

\begin{pf*}{Proof of uniqueness}
Let $v^1$ and $v^2$ be two bounded classical solutions to
the term structure equation, and define
\[
v(x,t)=v^1(x,T-t)-v^2(x,T-t).
\]
Then $v(x,t)$ is a bounded solution to
%
\begin{equation}
\label{eqn}
\cases{
v_t=\frac{1}{2}\sigma^2 v_{xx}+\beta v_x-xv,\vspace*{2pt}\cr
v(x,0)=0,\cr
v_t(0,t)=\beta(0,t)v_x(0,t).
}
\end{equation}
Now consider the function
\[
h(x,t)=(1+x)e^{Mt},
\]
where $M$ is a positive constant.
For $M$ large enough, depending on $\beta(0,t)$ and
the growth rate of $\beta$, $h$ is a super-solution to (\ref{eqn})
which tends to
infinity at spatial infinity. Thus,
according to the maximum principle, the function
$v$ is bounded above by $\varepsilon h$
and below by $-\varepsilon h$ for any $\varepsilon>0$. It follows that
$v\equiv0$,
which demonstrates uniqueness of bounded classical solutions
to the term structure equation.
\end{pf*}

\begin{pf*}{Proof of continuity}
To show that $u$ is continuous, denote by
$X^{x,t}$ the solution to (\ref{X}) with initial condition
$X^{x,t}(t)=x$. Let $(x,t)$ and $(y,r)$ be two points in
$[0,\infty)\times[0,T]$. Then, if
$r\leq t$, we have
%
\begin{eqnarray}
\label{uc}
\vert u(y,r)-u(x,t)\vert&\leq&
E\bigl[ e^{-\int_r^T X^{y,r}(s) \,ds}\vert g(X^{y,r}(T))-
g(X^{x,t}(T))\vert\bigr]\nonumber
\\
&& {}
+ E\bigl[g(X^{x,t}(T))\bigl\vert e^{-\int_r^T X^{y,r}(s)\, ds}
-e^{-\int_t^T X^{x,t}(s)\, ds} \bigr\vert\bigr]\nonumber
\\
&\leq& E[\vert g(X^{y,r}(T))-
g(X^{x,t}(T))\vert]
\\
\nonumber
&&{} +C\int_t^TE[\vert X^{y,r}(s)-X^{x,t}(s)\vert] \,ds
\\
\nonumber
&&{} +C\int_r^tE[X^{y,r}(s)]\, ds
\end{eqnarray}
for some constant $C$, where we have used that $g$ is bounded.
A similar expression can be derived if $r>t$.
It follows from Remark 1 in Section 8, Chapter 2 in \cite{GS} that
$X^{y,r}(t)\to x$ in $L^2$ as $(y,r)\to
(x,t)$. Therefore, from Theorem~2.1 in \cite{BMO}
we have
\[
E\Bigl[\sup_{t\leq s\leq T}\bigl(X^{y,r}(s)- X^{x,t}(s)\bigr)^2\Bigr]\to0
\]
as $(y,r)\to(x,t)$. (Theorem~2.1 in \cite{BMO}
also holds in the case of random starting points.)
Since $g$ is assumed continuous and bounded, all three terms
on the right-hand side of (\ref{uc}) tend to 0 as $(y,r)\to
(x,t)$. Thus $u$ is continuous on $[0,\infty)\times[0,T]$.
\end{pf*}

\begin{pf*}{Proof that $u\in C^{2,1}((0,\infty)\times[0,T))$ and
satisfies (\ref{tse})}
For a given point $(x,t)\in(0,\infty)\times[0,T)$,
let
\[
R=(x_1,x_2)\times[t_1,t_2)\subseteq(0,\infty)\times[0,T)
\]
be a rectangle which contains $(x,t)$, where $x_1>0$. Since $u$ is
continuous, it follows from standard parabolic theory, see \cite{Fr},
that there exists a unique solution
$U\in C^{2,1}(R)$ to the boundary value problem
\[
\cases{ U_t+\frac{1}{2}\sigma^2U_{xx}+\beta U_x-xU=0,&\quad
\mbox{in }$R$,\cr
U=u, &\quad\mbox{on }$\partial_p R$,
}
\]
where $\partial_p R=( [x_1,x_2]\times\{t_2\})\cup(\{x_1,x_2\}\times
[t_1,t_2])$ is the parabolic boundary of~$R$.
From It\^o's formula, the process
\[
Z(s)=e^{-\int_t^s X^{x,t}(r)\, dr}U(X^{x,t}(s),s)
\]
is a martingale on the time interval $[t,\tau_R]$, where
\[
\tau_R=\inf\{s\geq t\dvtx X^{x,t}(s)\notin R\}
\]
is the first exit time from the rectangle $R$.
Therefore,
\[
U(x,t)=E\bigl[e^{-\int_t^{\tau_R} X^{x,t}(r)\, dr}
u(X^{x,t}(\tau_R),\tau_R)\bigr]=u(x,t),
\]
where the second equality follows from the strong Markov
property. Consequently, $u\in C^{2,1}((0,\infty)\times[0,T))$.
Since $u\equiv U$ on $R$, we also see that $u$ satisfies
(\ref{tse}).
\end{pf*}

It remains to show that $u$ is continuously differentiable up
to the spatial boundary $x=0$, and that it satisfies the boundary
condition (\ref{bc}). This is
done in Sections~\ref{u_x}--\ref{sahlin}.

\section{Continuity of the first spatial derivative}
\label{u_x}

In this section we investigate regularity of the spatial
derivative $u_x$ at the boundary $x=0$. To do this
we study the stochastic representation of the terminal value
problem obtained by formally differentiating the term
structure equation. We show that this stochastic representation
indeed is the derivative of $u$ and that it is continuous.

Recall that $\alpha_x$ is assumed to be continuous on
$[0,\infty)\times[0,T]$, where $\alpha(x,t)= \frac{1}{2}\sigma
^2(x,t)$. Let the process $Y$ be
modeled by the stochastic differential equation
%
\begin{equation}
\label{Y}
dY(t)=(\alpha_x+\beta)(Y(t),t)\, dt + \sigma(Y(t),t)\, dW.
\end{equation}
Rather than specifying precise conditions under which (\ref{Y})
has a unique solution, we simply assume what we need.

\begin{assumption}
\label{ass}
The coefficients $\sigma$ and $\beta$ are
such that, path-wise, uniqueness holds for equation (\ref{Y}).
\end{assumption}

\begin{remark*}
 Note that Assumption~\ref{ass} holds for example if
$\alpha$ is
twice continuously differentiable in space, since then the drift
$\alpha_x+\beta$ is locally Lipschitz continuous.
Moreover,
if $\sigma$ and $\beta$ are time-independent,
then it follows from \cite{ABBP,BP}
and Section IX.3 in \cite{RY} that Assumption~\ref{ass}
automatically holds. Thus the
Cox--Ingersoll--Ross model (\ref{CIR}), the Dothan model
(\ref{D}), the Hull--White model (\ref{HW}) and the
model (\ref{exam}) all satisfy Assumption~\ref{ass}.

Also note that since $\alpha(0,t)=0$, we have
$\alpha_x(0,t)\geq0$.
Thus $Y$ remains nonnegative since it has the same volatility
as $X$ but a larger drift at $0$.
\end{remark*}

Next, define the function $v$ by
\begin{eqnarray}
\label{v}
v(x,t) &=& E\biggl[g^\prime(Y(T))
\exp\biggl\{\int_t^T \beta_x(Y(s),s)-Y(s) \,ds\biggr\}\biggr]\nonumber
\\[-8pt]\\[-8pt]
&& \nonumber
-E\biggl[\int_t^T
\exp\biggl\{\int_t^s \beta_x(Y(r),r)-Y(r) \,dr\biggr\}u(Y(s),s)\, ds\biggr],
\end{eqnarray}
where $Y$ is the solution to (\ref{Y}) with initial condition
$Y(t)=x$.

If the term structure equation (\ref{tse}) is formally
differentiated with respect to $x$, then the derivative
$u_x$ satisfies
\[
(u_x)_t+\alpha(u_x)_{xx}+(\alpha_x+\beta)(u_x)_x
+(\beta_x-x)u_x-u=0
\]
with terminal condition $u_x(x,T)=g^\prime(x)$. The function
$v$ defined in (\ref{v}) is the corresponding stochastic
representation. In Theorem~\ref{girsanov} below we show that $v$
indeed equals the spatial
derivative of $u$.

\begin{proposition}
\label{vcont}
The function
$v(x,t)$ is continuous on $[0,\infty)\times[0,T]$.
\end{proposition}

\begin{pf}
The result follows along the same lines as the continuity of $u$ above.
Indeed, let $(x_n,t_n)$
converge to $(x,t)$, where $t_n\leq t$,
and let $Y$ and $Y^n$ be defined by
\[
\cases{
dY(s)=(\alpha_x+\beta)(Y(s),s)\, ds + \sigma(Y(s),s) \,dW,\cr
Y(t)=x
}
\]
and
\[
\cases{
dY^n(s)=(\alpha_x+\beta)(Y^n(s),s)\, ds + \sigma(Y^n(s),s)\, dW,\cr
Y^n(t_n)=x_n,
}
\]
respectively. Also define
\[
I(s):= \exp\biggl\{ \int_t^s \beta_x(Y(u),u)-Y(u) \,du\biggr\}
\]
and
\[
I^n(s):= \exp\biggl\{ \int_{t_n}^s \beta_x(Y^n(u),u)-Y^n(u) \,du\biggr\}.
\]
Then
\begin{eqnarray*}
\vert v(x_n,t_n)-v(x,t)\vert&\leq&
E [\vert I^n(T) g^\prime(Y^n(T))-I(T) g^\prime(Y(T))\vert
]\\
&&{} + \int_t^T E[\vert I^n(s)u(Y^n(s),s)-I(s)u(Y(s),s)
\vert]\, ds\\
&&{} + \int_{t_n}^t E[ I^n(s) u(Y^n(s),s)]\, ds.
\end{eqnarray*}
The first term and the integrand in the second term are similar
to the type of terms treated when proving the continuity of
$u$. Moreover, the integrand of the third term is bounded.
Thus it follows from bounded convergence that $v$ is continuous.
\end{pf}

We also need a continuity result in the volatility parameter.
To formulate it, let
$\{\sigma^n(x,t)\}_{n=1}^\infty$ be a sequence of functions
satisfying Hypothesis~\ref{hyp} uniformly in $n$, that is, with the
same constant $C$ in the bound (\ref{bound}). Moreover,
assume that $\sigma^n(x,t)$ converges
to $\sigma(x,t)$ and $\alpha^n_x$ converges to $\alpha_x$ uniformly
on compacts as $n\to\infty$, where
$\alpha^n=\frac{1}{2}(\sigma^n)^2$.
Let $u^n$ and $v^n$ be defined as $u$ and $v$ but using the volatility
function $\sigma^n$ instead of $\sigma$. More explicitly,
\[
u^n(x,t)=E\bigl[e^{-\int_t^T X^n(s) \,ds}g(X^n(T))\bigr]
\]
and
\begin{eqnarray}
\label{vn}
v^n(x,t) &=& E\biggl[g^\prime(Y^n(T))
\exp\biggl\{\int_t^T \beta_x(Y^n(s),s)-Y^n(s)\, ds\biggr\}\biggr]\nonumber
\\[-8pt]\\[-8pt]
&& \nonumber{}
-E\biggl[\int_t^T \exp\biggl\{\int_t^s \beta_x(Y^n(r),r)-Y^n(r)\, dr\biggr\}
u^n(Y^n(s),s) \,ds\biggr],
\end{eqnarray}
where $X^n$ and $Y^n$ satisfy
\[
\cases{
dX^n(s)=\beta(X^n(s),s) \,ds + \sigma^n(X^n(s),s)\, dW(s),\cr
X^n(t)=x
}
\]
and
\[
\cases{
dY^n(s)=(\alpha^n_x+\beta)(Y^n(s),s) \,ds + \sigma^n(Y^n(s),s) \,dW(s),\cr
Y^n(t)=x,
}
\]
respectively.

\begin{proposition}
\label{parcont}
The functions $u$ and $v$ are continuous in the volatility parameter. More
precisely, $u^n(x,t)\to u(x,t)$ and
$v^n(x,t)\to v(x,t)$ as $n\to\infty$ for any fixed point
$(x,t)\in[0,\infty)\times[0,T]$.
\end{proposition}

\begin{pf}
It follows from Theorem~2.5 in \cite{BMO} that
\[
\lim_{n\to\infty}E\Bigl[\sup_{s\in[t,T]}\bigl(X(s)-X^n(s)\bigr)^2\Bigr]
=0.
\]
Therefore,
\begin{eqnarray*}
\vert u^n(x,t)-u(x,t)\vert&\leq&
E\bigl[\bigl\vert e^{-\int_t^TX^n(s) \,ds}-e^{-\int_t^TX(s)\, ds}\bigr\vert
g(X^n(T))\bigr]
\\
&&{} + E\bigl[e^{-\int_t^TX(s) \,ds}\vert g(X^n(T))-g(X(T))\vert\bigr]
\\
&\leq& C\int_t^TE[\vert X(s)-X^n(s)\vert] \,ds + E[\vert g(X^n(T))-g(X(T))\vert]
\\
&\to&0
\end{eqnarray*}
as $n\to\infty$. Thus $u$ is continuous in the volatility
function.

The continuity of $v$ in the volatility function is similar.
Indeed, let
\[
I(s):=\exp\biggl\{\int_t^s\beta_x(Y(r),r)-Y(r) \,dr\biggr\}
\]
and
\[
I^n(s):=\exp\biggl\{\int_t^s\beta_x(Y^n(r),r)-Y^n(r)
\,dr\biggr\}.
\]
Then
\begin{eqnarray*}
\vert v^n(x,t)-v(x,t)\vert&\leq&
E[\vert I^n(T) g^\prime(Y^n(T))-I(T)g^\prime(Y(T))\vert]
\\
&{}& +\int_t^T E[\vert I^n(s)u^n(Y^n(s),s)-I(s) u^n(Y(s),s)\vert]\, ds
\\
&&{} + \int_t^T E[I(s)\vert u^n(Y(s),s)-u(Y(s),s)\vert
] \,ds.
\end{eqnarray*}
The first term and the integrand of the third term
are similar to the terms
appearing in the first part of the proof. The integrand
of the second term can be dealt with using the fact that
each $u^n$ is Lipschitz continuous in $x$, uniformly in $n$ since the Lipschitz property is inherited
by the value function.
Thus all terms tend to zero as $n\to\infty$, so $v$ is
continuous in the volatility.
\end{pf}

\begin{theorem}
\label{girsanov}
We have $u_x(x,t)=v(x,t)$
on $[0,\infty)\times[0,T]$. Consequently, $u_x$ is continuous
on $[0,\infty)\times[0,T]$.
\end{theorem}

\begin{pf}
It suffices to prove $u_x(x,0)=v(x,0)$. We first assume
that $\sigma$ is continuously differentiable in $x$ with
a bounded derivative. It then follows
from Section~5.5 in \cite{F} or Section 8 in \cite{GS} that the derivative
\[
\xi(t):=\frac{\partial X(t)}{\partial x}
\]
of $X(t)=X^{x,0}(t)$ with respect to the initial
point $x$ exists and is continuous, and it satisfies
\[
\cases{
d\xi(t)= \xi(t)\beta_x(X(t),t) \,dt +
\xi(t)\sigma_x(X(t),t) \,dW(t),\cr
\xi(0)=1.
}
\]
Moreover,
%
\begin{eqnarray}
\label{ux}
u_x(x,0) &=& E\biggl[g^\prime(X(T))\xi(T)
\exp\biggl\{-\int_0^T X(s) \,ds\biggr\}\biggr]\nonumber
\\
&&
{}-E\biggl[ g(X(T)) \exp\biggl\{-\int_0^T X(s) \,ds\biggr\} \int_0^T\xi(s)\, ds\biggr]
\\
&=:& \nonumber I_1-I_2.
\end{eqnarray}
We claim that
$I_i=J_i$, $i=1,2$, where
\[
J_1 = E\biggl[g^\prime(Y(T))
\exp\biggl\{\int_0^T \beta_x(Y(s),s)-Y(s) \,ds\biggr\}\biggr]
\]
and
\[
J_2 =E\biggl[\int_0^T
\exp\biggl\{\int_0^s \beta_x(Y(r),r)- Y(r)\, dr \biggr\}u(Y(s),s)\, ds\biggr],
\]
compare (\ref{v}) above. Here $Y$ is defined as in (\ref{Y})
with initial condition $Y(0)=x$.

To show that $I_1=J_1$, define a new measure $Q$ on $\mathcal F_T$ by
$dQ=M(T)\, dP$, where the process $M$ is defined by
%
\begin{equation}
\label{M}
M(t)=\xi(t)\exp\biggl\{-\int_0^t \beta_x(Y(s)) \,ds\biggr\}.
\end{equation}
By It\^o's formula,
\[
dM(t)=M(t)\sigma_x(X(t))\, dW(t),
\]
so $M$ is a martingale since $\sigma_x$ is bounded.
In particular, $E[M(T)]=1$, so
$Q$ is a probability measure. From Girsanov's theorem it follows
that
\[
\tilde W(t)=W(t)-\int_0^t\sigma_x(X(s)) \,ds
\]
is a $Q$-Brownian motion, and
\[
dX=(\sigma\sigma_x+\beta)(X(t),t)\, dt +
\sigma(X(t),t)\, d\tilde W.
\]
Here $\sigma\sigma_x=\alpha_x$, so by weak uniqueness, the $Q$-law
of $X$ is the same as the law
of $Y$ under $P$. Consequently,
\begin{eqnarray*}
I_1 &=& E\biggl[g^\prime(X(T))\xi^x(T)\exp\biggl\{-\int_0^TX(s)\, ds
\biggr\}\biggr]
\\
&=& E^Q\biggl[g^\prime(X(T))\exp\biggl\{\int_0^T \beta_x(X(s),s)-X(s) \,ds\biggr\}\biggr]
=J_1.
\end{eqnarray*}

To prove $I_2=J_2$, note that
\begin{eqnarray*}
I_2 &=& E\biggl[ g(X(T)) \exp\biggl\{-\int_0^T X(s)\, ds\biggr\} \int_0^T\xi(s)\, ds\biggr]
\\
&=& \int_0^T E\biggl[ \exp\biggl\{-\int_0^s X(r)\, dr\biggr\}\xi(s)
\\
&&\hspace*{28pt}  {}\times E\biggl[g(X(T))
\exp\biggl\{-\int_s^T X(r) \,dr\biggr\}\Big\vert\mathcal F_s\biggr]\biggr] \,ds
\\
&=& \int_0^T E\biggl[ \exp\biggl\{-\int_0^s X(r)\, dr\biggr\} \xi(s)u(X_s,s)\biggr]\, ds
\end{eqnarray*}
by the Markov property.
Define a new measure $Q=Q_s$ on $\mathcal F_s$ by
\[
dQ=M(s) \,dP,
\]
where $M$ is defined as in (\ref{M}). Girsanov's theorem
yields
\begin{eqnarray*}
&&  E\biggl[\exp\biggl\{-\int_0^s X(r)\, dr\biggr\}
\xi(s)u(X_s,s)\biggr]
\\
&&\qquad= E^Q\biggl[\exp\biggl\{\int_0^s\beta_x(X(r),r)-X(r)\, dr\biggr\} u(X_s,s)\biggr]
\\
&&\qquad = E\biggl[ \exp\biggl\{\int_0^s \beta_x(Y(r),r)-Y(r)\, dr\biggr\}u(Y_s,s)\biggr].
\end{eqnarray*}
Consequently, $I_2=J_2$, which finishes the proof
in the case of continuously differentiable $\sigma$.

The general case follows by approximation.
Let $\sigma^n$, $u^n$ and $v^n$ be as described before
Proposition~\ref{parcont}, with each $\sigma^n$ being
continuously differentiable in $x$ with bounded derivative.
From above, we then know that $v^n(x,t)=u^n_x(x,t)$ at all points.
Moreover, by Proposition~\ref{parcont}, $v^n(x,t)\to v(x,t)$
point-wise as $n\to\infty$.

On the other hand, since $u^n$ converges to $u$ point-wise and
is uniformly bounded, it
follows from standard parabolic theory that also $u_x^n$
converges to $u_x$ point-wise for all points $(x,t)$ with $x>0$.
Consequently, $v=u_x$ on $(0,\infty)\times[0,T]$. Since $v$ is
continuous on $[0,\infty)\times[0,T]$ by Proposition~\ref{vcont},
it is easy to check that $u_x(0,t)$ exists and that
we have $v=u_x$
everywhere on $[0,\infty)\times[0,T]$. The continuity of $u_x$
thus follows.
\end{pf}

\section{An estimate of the second spatial derivative}
\label{u_xx}

Since the function $v$ defined in (\ref{v}) is continuous,
it follows that [by a similar argument as in the proof that
$u$ satisfies (\ref{tse})] it indeed solves the
differentiated equation
\[
v_t=\alpha v_{xx}+(\alpha_x+\beta) v_x+(\beta_x-x)v-u
\]
on $(0,\infty)\times[0,T)$.
In this section we use interior estimates to show that $\alpha v_{x}\to
0$ as $x\to0$. Since $v=u_x$ by Theorem~\ref{girsanov},
this shows that the term $\alpha u_{xx}$ in (\ref{tse})
approaches zero close to the boundary.

\begin{proposition}
\label{uxx}
The function $v=u_x$ satisfies
\[
\lim_{(x,t)\to(0,t_0)}\alpha(x,t) v_x(x,t)=0
\]
for any $t_0$.
Consequently, $\lim_{(x,t)\to(0,t_0)}\alpha(x,t) u_{xx}(x,t)=0$.
\end{proposition}

\begin{pf}
Let $\{(x_n,t_n)\}_{n=1}^\infty\subseteq(0,\infty)\times[0,T)$
be a sequence of points converging to $(0,t_0)$, where
$t_0\in[0,T)$. Define new coordinates $(y,s)$
by letting $y=kx$ and $s=k(t-t_0)$, where $k$ is specified
more precisely below. Then the function $w$ defined by
\[
w(y,s)=v(x,t)
\]
satisfies
%
\begin{equation}
\label{weq}
w_s=\tilde\alpha w_{yy} + \tilde\beta
w_y+\gamma w+ h,
\end{equation}
where
\begin{eqnarray*}
\tilde\alpha(y,s)&=&\alpha\biggl(\frac{y}{k},t_0+\frac{s}{k}\biggr)k,
\\
\tilde\beta(y,s)&=&(\alpha_x+\beta)\biggl(\frac{y}{k},t_0+\frac{s}{k}\biggr),
\\
\gamma(y,s)&=&\frac{1}{k}\beta_x\biggl(\frac{y}{k},t_0+\frac{s}{k}\biggr)-\frac{y}{k^2}
\end{eqnarray*}
and
\[
h(y,s)=-\frac{1}{k}u\biggl(\frac{y}{k},t_0+\frac{s}{k}\biggr).
\]
Now consider a region $\mathcal R=\mathcal R^n$ which contains
the point $(x_n,t_n)$, and such that
%
\begin{equation}
\label{R}
1\leq\alpha(x,t)k\leq2
\end{equation}
in $\mathcal R$. Since $\alpha_x(x,t)$ is continuous up to the
boundary, the region $\mathcal R$ in $(y,s)$-coordinates does not
collapse as $n\to\infty$, but it can rather be chosen to
consist of a rectangle of
fixed size; the location of the rectangle is not necessarily
fixed though. In this rectangle, the coefficients of the equation
(\ref{weq})
satisfy
\begin{eqnarray*}
1&\leq&\tilde\alpha(y,s)\leq2,
\\
\vert\tilde\beta(y,s)\vert&\leq& C,
\\
\vert\gamma(y,s)\vert&\leq& C
\end{eqnarray*}
and
\[
\vert h(y,s)\vert\leq C/k
\]
for some constant $C$ which is independent of $n$.
Since $w(y,s)=v(x,t)$ we have that $w$ converges to
the constant $v(0,t_0) = u_x(0,t_0)$ uniformly on
$\mathcal R$ as $n\to\infty$. By interior Schauder estimates,
$w_y$ tends to 0 as $n\to\infty$. Since
\[
\alpha(x,t)v_x(x,t)=\tilde\alpha(y,s) w_y(y,s),
\]
and since $\tilde\alpha(y,s)$ is bounded on $\mathcal R$,
the conclusion follows.
\end{pf}

\section{The time derivative at the boundary}
\label{sahlin}

It follows from Proposition~\ref{uxx} and (\ref{tse}) that
%
\begin{equation}
\label{lim}
\lim_{(x,t)\to(0,t_0)} u_t(x,t)+\beta(0,t_0)u_x(0,t_0)=0
\end{equation}
for any $t_0\in[0,T)$. In this section we show that the boundary
condition (\ref{bc}) also holds \textit{at} the boundary, that is, not
merely in the limit.

\begin{proposition}
\label{atzero}
The function $u_t(x,t)+\beta(x,t)u_x(x,t)$ defines
a continuous function on $[0,\infty)\times[0,T)$. Moreover,
it vanishes for $x=0$.
\end{proposition}

\begin{remark*}
 Note that Proposition~\ref{atzero} finishes the
proof of
Theorem~\ref{main}.
\end{remark*}

\begin{pf*}{Proof of Proposition~\ref{atzero}}  In view of (\ref{lim})
above, it suffices to show that
$u_t$ exists at the boundary and that it equals $-\beta u_x$.
To do this, fix a point on the boundary with coordinates
$(0,t_0)$. For notational simplicity we assume that
$t_0=0$.
The time (left) derivative $u_t$ at the boundary is defined by
%
\begin{equation}
\label{limits}
u_t(0,0)=\lim_{k\to\infty}k\biggl( u(0,0)-u\biggl(0,-\frac{1}{k}\biggr)\biggr),
\end{equation}
provided the limit exists. To determine $u_t(0,0)$,
we let $X^k$ be defined by
\[
\cases{
dX^k=\beta(X^k(t),t) \,dt +\sigma(X^k(t),t)\, dW,\cr
X^k(-1/k)=0.
}
\]
However, instead of considering the process $X$ with
different starting times, we perform a change of variables
so that the starting time is independent of $k$. We thus
introduce the process
$Y^k(s)$ by
\[
Y^k(s)= kX^k\biggl(\frac{s}{k}\biggr).
\]
With respect to the time variable $s$, the dynamics
of $Y^k$ has the form
%
\begin{equation}
\label{Yk}
\cases{
dY^k(s)=\beta\biggl(\dfrac{1}{k}Y^k(s),\dfrac{s}{k}\biggr) \,ds+ \sqrt{k\sigma
^2\biggl(\dfrac{1}{k}Y^k(s),\dfrac{s}{k}\biggr)} \,dW^k,\cr
Y^k(-1)=0,
}
\end{equation}
where $W^k(s)$ denotes some Brownian motion.
By the Markov property,
\begin{eqnarray*}
u(0,-1/k) &=& E\bigl[e^{-\int_{-1/k}^0{X^k(s)} \,ds}
u(X^k(0),0)\bigr]
\\
&=& E\biggl[e^{-\int_{-1}^0{{(1/k^2)}Y^k(s)} \,ds}
u\biggl(\frac{1}{k}Y^k(0),0\biggr)\biggr].
\end{eqnarray*}
Hence,
\begin{eqnarray*}
u_t(0,0) &=&\lim_{k\to\infty} kE_{0,-1} \biggl[u(0,0)-e^{-\int
_{-1}^0{{(1/k^2)}Y^k(s)}\, ds}
u\biggl(\frac{1}{k}Y^k(0),0\biggr)\biggr]
\\
&=&\lim_{k\to\infty} E_{0,-1} \biggl[k \biggl(u(0,0)-u\biggl(\frac{1}{k}Y^k(0),0\biggr)\biggr)\biggr],
\end{eqnarray*}
where the second equality follows using the inequality $e^{-x}-1\ge-x$ since
\begin{eqnarray*}
E_{0,-1} \biggl[ku\biggl(\frac{1}{k}Y^k(0),0\biggr)\bigl\vert e^{-\int_{-1}^0{
{(1/k^2)}Y^k(s)}\, ds}-1\bigr\vert\biggr]
\leq C\frac{1}{k} E_{0,-1}\biggl[\int_{-1}^0Y^k(s) \,ds\biggr]
\to 0
\end{eqnarray*}
as $k\to\infty$. Now, define the process $Y$ by
\[
\cases{
dY=\beta(0,0) \,ds + \sqrt{2\alpha_x(0,0)Y} \,dW,\cr
Y(-1)=0,
}
\]
and redefine $Y^k$ as in (\ref{Yk}) above but using the same
Brownian motion $W$  (this does not change the law of $Y^k$).
Since
\[
\beta^k(y,s):=\beta\biggl(\frac{y}{k},\frac{s}{k}\biggr)\to
\beta(0,0)
\]
and
\[
\sigma^k(y,s):=\sqrt{k\sigma^2\biggl(\frac{y}{k},\frac{s}{k}\biggr)}
\to\sqrt{2\alpha_x(0,0)y}
\]
uniformly on compacts as $ k\to\infty$
(here we used the assumption that\vspace*{1pt} $\alpha$ is
continuously differentiable in space),
it follows from \cite{BMO} that $Y^k(0)\to Y(0)$ in $L^2$ as
$k\to\infty$. From Theorem~\ref{girsanov} above we know
that $u$ is differentiable in $x$, so
$k(u(0,0)-u(\frac{y}{k},0))$ converges to $-u_x(0,0)y$. By dominated
convergence, we have
\[
kE \biggl[u(0,0)-u\biggl(\frac{1}{k}Y(0),0\biggr)\biggr]\to
-u_x(0,0)E[Y(0)]=-\beta(0,0)u_x(0,0)
\]
as $k\to\infty$. Moreover, the Lipschitz property of $u$ yields
that
\[
kE\biggl[u\biggl(\frac{1}{k}Y(0),0\biggr)-
u\biggl(\frac{1}{k}Y^k(0),0\biggr)\biggr] \leq
C E[ \vert Y(0)-Y^k(0)\vert]\to0
\]
as $k\to\infty$. It follows that
\begin{eqnarray*}
u_t(0,0)+u_x(0,0)\beta(0,0)=0 .
\end{eqnarray*}
As $t_0=0$ was chosen only for notational convenience,
we have that
\begin{eqnarray*}
u_t(0,t)+\beta(0,t)u_x(0,t)=0
\end{eqnarray*}
for any $t$.
To be precise, we have shown the result above only for the left $t$-derivative.
However, this left $t$-derivative is continuous by the equation above,
so it follows from a simple calculus lemma
that in fact $u$ is differentiable in time, thus finishing our proof.
\end{pf*}

\section{Models allowing negative interest rates}
\label{realline}

For models in which the short rate can fall below zero with
positive probability, the connection between
the option price, given by a risk-neutral
expected value, and the term structure equation is more
straightforward than for models with nonnegative rates.
Nevertheless, we have not been able to find a precise reference
for this case, so for completeness we provide such a result
in this section.

The assumptions needed on the drift and volatility are presented
below, where $x^+=\max(x,0)$. These assumptions now replace
those of Hypothesis~\ref{hyp}.
To our knowledge, all models used in practice that allow negative
interest rates satisfy these requirements.
For example, the
Vasicek model, in which
\[
dX(t)=\bigl(a-bX(t)\bigr) \,dt + \sigma\, dW,
\]
where $a$, $b$ and $\sigma$ are positive constants, is covered.

\begin{hypothesis}
\label{whole}
The drift $\beta\in C(\mathbb{R}\times[0,T])$ and the volatility
$\sigma\in
C(\mathbb{R}\times[0,T])$ are both
Lipschitz continuous in $x$. Moreover,
\[
0<\vert\sigma(x,t)\vert\leq C(1+x^+)
\]
and
%
\begin{equation}
\label{C}
\vert\beta(x,t)\vert\leq C(1+\vert x\vert)
\end{equation}
for some positive constant $C$.
\end{hypothesis}

\begin{remark*}
 The assumption that $\sigma$ is strictly positive is
not a strong assumption. Indeed, let us for simplicity
consider a time-homogeneous model $dX(t)=\beta(X(t)) \,dt + \sigma
(X(t)) \,dW(t)$ with $\sigma(a)=0$ for some $a\in\mathbb{R}$.
If $X(0)\geq a$ and $\beta(a)\geq0$, then the process $X$ cannot
take values smaller than $a$, so we are essentially in the situation
handled in Sections~\ref{halfline}--\ref{sahlin} (but with the point 0
replaced with $a$). If $b(a)<0$, then $X$
can take values below $a$, but if this happens then
the process will stay below $a$ forever, and we are then
again back in the previous situation.

The bound on the volatility for negative rates guarantees
that bond prices are finite. For models in which $\sigma$ grows
faster than $\sqrt{\vert x\vert}$ for negative rates, bond
prices can be infinite; compare Theorem~4.1 in \cite{Y}.
\end{remark*}

For a given continuous payoff function $g\dvtx\mathbb{R}\to[0,\infty)$,
define the corresponding option price $u\dvtx\mathbb{R}\times[0,T]$ by
\[
u(x,t)=E_{x,t}\biggl[\exp\biggl\{
-\int_t^T X(s) \,ds\biggr\}g(X(T))\biggr],
\]
where
\[
\cases{
dX(s)=\beta(X(s),s)\, ds + \sigma(X(s),s) \,dW(s),\cr
X(t)=x.
}
\]
We require that the payoff function is bounded for positive
interest rates and of, at most, exponential growth for
negative rates, that is,
%
\begin{equation}
\label{g}
0\leq g(x)\leq K \max\{1,e^{-Kx}\}
\end{equation}
for some positive constant $K$. The corresponding term structure
equation is given by
\begin{eqnarray*}
u_t(x,t)+\tfrac{1}{2}\sigma^2(x,t) u_{xx}(x,t)+\beta(x,t) u_x(x,t)=xu(x,t)
\end{eqnarray*}
on $\mathbb{R}\times[0,T)$, with terminal condition
$u(x,T)=g(x)$. By a {\it classical solution}
to this equation we mean a solution which is continuous up to the boundary
$t=T$ and with all derivatives appearing in the equation being
continuous functions on the set $t<T$.
The following is our main result
about the term structure equation on the whole real line.

\begin{theorem}
\label{palme}
Assume Hypothesis $ \ref{whole}$\ and the bound (\ref{g}). Then
the option price $u(x,t)$ satisfies
%
\begin{equation}
\label{ubound}
u(x,t)\leq K^\prime\max\{1,e^{-K^\prime x}\}
\end{equation}
for some constant $K^\prime$. Moreover,
$u(x,t)$ is the unique classical solution to
the term structure equation satisfying this growth assumption for some
constant $K^\prime$.
\end{theorem}

\begin{remark*}
 Note that the bound (\ref{g}) is natural for models
on the whole real line. In fact, even if $g$ was
bounded, the option price $u$ would be of exponential
growth for negative rates. Also note that, for example,
call options on a bond are covered by Theorem~\ref{palme}.
Indeed, the payoff of a bond call
option with maturity $T_1$ is given by $g(x)=(u(x,T_1)-K)^+$,
where $u$ is the price of a bond maturing at $T_2>T_1$.
Since $u$ satisfies (\ref{ubound}) by Theorem~\ref{palme},
the payoff $g$ satisfies (\ref{g}).
Bond put options are trivially covered since they have
bounded payoff functions.
\end{remark*}

\begin{pf*}{Proof of Theorem~\ref{palme}} The bound (\ref{ubound})
follows from Corollary~3.3
in \cite{ET}. To prove uniqueness of solutions, assume that $v$ is a solution
to the term structure equation with boundary value $g=0$
such that $\vert v\vert\leq K^\prime\max\{1,e^{-K^\prime x}\}$
for some constant $K^\prime$. Let
\[
h(x,t)= e^{M(T-t)}\bigl(e^{-f(t)x}+x\bigr)
\]
for some large constant $M$. Here
\[
f(t)=\frac{e^{C(T-t)}-1}{C}+ Ke^{C(T-t)},
\]
where $C$ is the constant appearing in (\ref{C}) and
$K>K^\prime$. Then the set
\begin{eqnarray*}
\label{h1}
\{(x,t)\in\mathbb{R}\times[0,T]\dvtx\varepsilon h(x,t)< v(x,t) \}
\end{eqnarray*}
is bounded, and
\[
h_t+\tfrac{1}{2}\sigma^2h_{xx}+\beta h_x-xh<0
\]
at all points provided $M$ is chosen large enough.
Standard methods used to prove the
maximum principle yield that $v=0$ at all points. Thus we have
uniqueness of solutions to the term structure equation
in the class of functions satisfying (\ref{ubound}).

To show that $u$ is a classical solution to the term
structure equation, we carry out an approximation argument. We consider
the term structure equation but with the discount
factor $x$ replaced by bounded functions that agree with
$x$ inside some large compact sets. We also
replace the payoff function
$g$ with functions of, at most, polynomial growth that agree with $g$
on large compact sets. The corresponding equation
is in the standard class and its stochastic solution is known to be
continuous and hence is a classical solution. Now let the functions
approximating the
discount factor grow up to $x$ for large positive $x$ and then decrease
down to $x$ for large negative $x$. By the monotone convergence theorem,
the corresponding
stochastic solutions converge to the stochastic solution above denoted
$u$. Interior Schauder estimates yield interior regularity of the
limiting solution and
continuity at
the boundary is established using the maximum principle.
\end{pf*}

%

\printaddresses

\end{document}